\documentclass[11pt,reqno,fleqn]{article}
\usepackage[T1]{fontenc}
\usepackage{lmodern}
\usepackage[scaled]{}
\usepackage[ngerman,french,english]{babel}
\usepackage[utf8]{inputenc}

\usepackage{epsfig}
\usepackage{graphicx}
\usepackage{wrapfig}
\usepackage{amsmath}

\usepackage{amssymb}
\usepackage{amscd}
\usepackage{amsthm}
\usepackage{tikz}
\usetikzlibrary{arrows}
\newcommand{\bq}{\begin{equation}}
\newcommand{\eq}{\end{equation}}
\newcommand{\ba}{\begin{align}}
\newcommand{\ea}{\end{align}}
\newcommand{\be}{\begin{eqnarray}}
\newcommand{\ee}{\end{eqnarray}}
\newcommand{\bi}{\begin{itemize}}
\newcommand{\ei}{\end{itemize}}
\newcommand{\bd}{\begin{definition}}
\newcommand{\ed}{\end{definition}}
\newcommand{\bt}{\begin{theorem}}
\newcommand{\et}{\end{theorem}}
\newcommand{\bc}{\begin{corollary}}
\newcommand{\ec}{\end{corollary}}
\newcommand{\bcn}{\begin{conjecture}}
\newcommand{\ecn}{\end{conjecture}}
\newcommand{\br}{\begin{remark}}
\newcommand{\er}{\end{remark}}
\newcommand{\ce}{\begin{eqnarray*}}
\newcommand{\de}{\end{eqnarray*}}
\newcommand{\bpf}{\begin{proof}}
\newcommand{\epf}{\end{proof}}
\newcommand{\bl}{\begin{lemma}}
\newcommand{\el}{\end{lemma}}
\newtheorem{theorem}{Theorem}[section]

\newtheorem{conjecture}[theorem]{Conjecture}
\newtheorem{lemma}[theorem]{Lemma}
\newtheorem*{remark}{Remark}
\newtheorem{definition}[theorem]{Definition}
\newtheorem{proposition}[theorem]{Proposition}
\newtheorem{Examples}[theorem]{Examples}
\newtheorem{corollary}[theorem]{Corollary}
\numberwithin{equation}{section}

\def\b{\beta}
\def\d{\delta}

\def\[{{\Big[}}
\def\]{{\Big]}}
\def\<{{\langle}}
\def\>{{\rangle}}
\def\s{\sigma}
\def\({{\Big(}}
\def\){{\Big)}}

\def\bt{\begin{theorem}}
\def\et{\end{theorem}}
\def\bl{\begin{lemma}}
\def\el{\end{lemma}}
\def\br{\begin{remark}}
\def\er{\end{remark}}
\def\bx{\begin{Examples}}
\def\ex{\end{Examples}}
\def\bd{\begin{definition}}
\def\ed{\end{definition}}
\def\bp{\begin{proposition}}
\def\ep{\end{proposition}}
\def\bc{\begin{corollary}}
\def\ec{\end{corollary}}

\def\cB{{\mathcal B}}

\def\cN{{\mathcal N}}

\def\cS{{\mathcal S}}

\def\mN{{\mathbb N}}

\def\mX{{\mathbb X}}

\def\mZ{{\mathbb Z}}

\usepackage{graphicx}

\begin{document}

\section* {\center{A Family of Totally Rank One Two-sided Shift Maps}}
\begin{center}
\textsc
{
Yue Wu,\footnote{Schlumberger WesternGeco Geosolution, Houston, Texas, USA}
Dongmei Li,\footnote{Corresponding Author; Harbin University of Science and Technology, School of Applied Science, Harbin, Heilongjiang, China}
Yunjian Wang,\footnote{Schlumgerger SWT Technology Center, Houton, Texas, USA}
Diquan Li,\footnote{Central South University, School of Geosciences and Info-physics,University, Changsha, Hunan, China}
}
\end{center}

\selectlanguage{english}
\begin{abstract}
'Generalized del Junco-Rudolph's map', a sub-family of generalized Chacon's map (\cite{FER1}), is introduced. A skew product related to the structure of the Generalized del Junco-Rudolph's map is introduced. A Relative Prime Relation, $h_{k+1}\equiv 1 \mod q$ is verified, based on the proposition of iterations of this skew product. We say a measure-preserving transformation $T$ is totally rank one if $T^{n}, n\in \mN$ is rank one. In this paper, we show that of every Generalized del Junco-Rudolph's map is totally rank one.
\end{abstract}
\section* {}
Since Boltzmann introduced the term of 'Ergodic' in his work on statistical mechanics, the ergodic theory has been studied and developed actively regarding modelized physical systems with dynamic nature. From ergodicity, mixing properties, isomorphism of systems, to simplicity, and many other properties, we seriously desire concrete examples to better understand and study these properties with related classifications.
The cutting-and-stacking method provides effective measure theoretical approximations of transformations, and at the same time, it provides a fruitful generating machine for typical examples and counter examples corresponding to symbolic shifting maps. Extending from the original construction by R.V.Chacon (\cite{CHA}), A.A. Prikhod'ko and V.V. Ryzhikov (Construction 7 in \cite{PRI1}) further described how to construct a transformation on $[0,1]$ using the cutting-and-stacking method. This kind of constructional approach may help us understand not only the symbolic representing of a dynamic system but also its powers. G.R. Goodson (\cite{GOO1}) studied the condition for an ergodic automorphism $T$ to be conjugate to its composition square $T^{2}$. What's more del Junco (\cite{JUN2}) showed that the property of being conjugate to its square, is a non-generic property of automorphisms. Thus therefore, it is not trivial at all to study the powers of ergodic automorphisms.\\

    We say a finite measure-preserving system $(\mX, \b, T, \mu)$ is rank one if a series of Rokhlin tower with base subset $F$ approximates any partition $P$ of $\mX$. Rank one system can be defined in multiple literatures (\cite{FER1}).
Rank one implies simple spectrum (\cite{BAX},\cite{GLA}), rigidity implies
singular spectrum (\cite{GLA}). The cutting and stacking structure also
gives us a way to study a rank one system by the symbolic
methods(\cite{CHA}, \cite{JUN1},etc.) All these make it
interesting to learn whether a transformation is rank one. Veech \cite{VEE3} showed that, measure theoretically, almost all interval exchange transformations are rank one. A famous rank one transformation, the Chacon's map coded as $B_{k+1}=B_{k}B_{k}1B_{k}$, and its extension are also studied. In this
chapter we will study the family of generalized del Junco-Rudolph's maps, and extend the rank one property of this family of maps (\cite{JUN1}), to their powers of natural numbers. We admit the notion of totally rank one for this property (Definition \ref{Total}), and make the main conclusion as Theorem \ref{main}. \\

\section{Generalized Del Junco-Rudolph's maps}
In \cite{JUN1}, an example described by symbolic recursion
is given. We extend this map, notated by Del Junco-Rudolph's map here, to a family of symbolic maps-generalized Del Junco-Rudolph's maps. It is a sub-family of the family of generalized Chacon's map \cite{FER1} We will describe the construction and some observations of the generalized Del Junco-Rudolph's maps in this section.\\
\\
The recursion formula determines the language of the system, thus
also the phase space of the system, which is a subset of $\mZ$
product space of the set of alphabet ${0,1}$. That is:\\
$$\begin{array}{l}
 B_{0}=0 \\
 B_{1}=(B_{0})^{a}1(B_{0})^{b-a}, (1\leq a \leq b-1, b\geq 2)\\
\cdots\\
B_{k+1}=(B_{k})^{a\cdot b^{k}}1(B_{0})^{(b-a)\cdot b^{k}}
\end{array}$$\\
Let $X\subset \{0,1\}^{\mZ}$ be defined by $X=\{x=\cdots
x_{-1}x_{0}x_{1}\cdots x_{k}\cdots |\,\mbox{for any }m\in \mZ, l\in
\mN\}$,
 $x_{m}x_{m+1}\cdots x_{m+l-1}$ is a
consecutive sequence of some $B_{k}$, $k\in \mN$. Let $T$ be the
left-shifting transformation on $X$, so $T$ is a two-sided
shift map.\\

\bd If  a sequence
$\overline{a} $ appears consecutively in some $B_{k}$, we say
$\overline{a}$ is
a valid sequence or a word in $\mX$ (language) of length $l$.
\ed

Next some combinatorical fact about $T$ will be shown, some
further concerns will be discussed in Chapter4 for the special
purpose. The following facts are easily seen or computed. Some of
them have been pointed out in [del Junco, Rudolpf
1].\\

\bl Let $T$ be the symbolic transformation defined above, i.e. the generalized Del.Junco-Rudolph map, we have the follow propositions:
    \begin{enumerate}
    \item $h_{k+1}=b^{k+1}h_{k}+1$, and
$b^{\frac{k(k+1)}{2}}<h_{k}<b^{\frac{k(k+1)}{2}+1}$.

    \item except for the two given copies of
$B_{n}$, there is no other consecutive sequence of letters which
is
identically $B_{n}$.

    \item suppose
$x|^{m+2h_{k}}_{m}=B_{k}1B_{k}$, then
$x|^{m-1}_{m-h_{k}}=B_{k}=x|^{m+3h_{k}}_{m+2h_{k}+1}$, for any $k>1$.\\
    \end{enumerate}
\el

\begin{proof}
(iii) Since $B_{k}1B_{k}$ is in some $B_{j}$, $j>k$, both the left
and the right word to $B_{k}1B_{k}$ of length $h_{k}$ are $B_{k}$.
\end{proof}
Topologically,we have the following conclusion:\\

\bl The system $(\mX,
T)$ is minimal.
\el
\begin{proof}
This is true based on the coding construction of the system.
\end{proof}
Measure theoretically a $T$-invariant measure of $\mX$ may be
determined by the measure of each cylinder, which is the
asymptotic density of the cylinder name in the language of the
symbolic structure. Del Junco Rudolpfh's approximation would be
introduced for this purpose, which asures the coincidence of the
measure defined and the asymptotic
density by the the uniform ergodic theorem.
\br We say $A$ is a cylinder set with name $\overline\b$ ($\left| \overline\b\right|=l$) at position $j$, if $A=\{x|x_{j}\cdots x_{j+l-1}=\overline \b\}$.
\er
Let $d_{k}(\overline{a})$ be the density in $B_{k}^{\mZ}=\cdots
B_{k}B_{k}\cdots$ of the occurrences of a word
$\overline{\b}=a_{0}\cdots a_{l-1}$. Then for $k$ large enough:
\be
\left|
d_{k+1}(\overline{\b})-d_{k}(\overline{\b})\right|<2l/h_{k+1}
\ee
This inequality is true because the only difference is taking
place at the higher
spacer (the spacer between two $B_{k}$).

\bd Let
$\cS_{k}=\{x:x|^{h_{k}}_{-h_{k}}=B_{k}1B_{k}\}$.
\ed

Next the measure $\mu$ is introduced based on the notion of $d_{k}(\overline\b)$:

\be
d(\overline{\b})=\lim_{k\rightarrow \infty}d_{k}(\overline{\b})
\ee
If $A$ is a cylinder set with name $\overline{\b}$, we define
\be
\mu (A)=d(\overline{\b})
\ee and extend
this countably additive measure to a shift invariant Borel measure
$\mu$ on $(\mX ,\cB, T)$, where $\cB$ is the $\s$-
algebra generated by the set of all the cylinder sets.\\

\bl
$\mu$ is the unique
invariant measure of $(\mX , \cB ,T)$
up to a multiple.
\el

\begin{proof}
Each point $x\in \mX$ is a generic point of $(\mX ,\mu ,\cB ,T)$.
\end{proof}
$\;$The measure of $\cS_{k}$ is bounded by:\\

\bp \label{ineq1}
$\frac{\displaystyle 1}{\displaystyle h_{k+1}}<\mu(\cS_{k})<\frac{\displaystyle b}{\displaystyle h_{k+1}}$
\ep
\begin{proof}
Occurrences of $\cS_{k}$ on the orbit of any $x\in\mX$ on the
orbit of any $x\in\mX$ are separated by at most $h_{k+1}$ and at
least $\frac{h_{k+1}+1}{b}$.
\end{proof}
\br
The above statement is labeled as a
proposition, since it is crucial for the future evaluation.
\er

The rigidity of $(\mX ,\cB ,\mu ,T )$ is easily seen by looking
at $T^{h_{k}}(C) \Delta C$ for any cylinder $C$, typically just consider
$B_{\Delta}^{n,k}=T^{h_{n}}(B_{k}) \Delta B_{k}, n\geq k$. The notion of coding distance provides us another way to
understand the rigidity of $(\mX ,\cB ,\mu ,T)$ and to evaluate
the approximate speed of the power of $T$ tending to indentity
along
the sequence $\{h_{k}\}$.\\

Let $\xi^{N,k}=B_{k}^{N}$ and set
\begin{align}
\d^{(N)}_{k}(T^{t}) & =\frac{1}{Nh_{k}}\sum\limits_{i=1}^{Nh_{k}}d^{(N,k)}_{code}(i,t)
\end{align}
where
\be
d^{(N,k)}_{code}(i,t)=\left\{\begin{array}{ll}
\left|\xi ^{N,k}(i+t)-\xi ^{N,k}(i)\right| & 1\leq i\leq Nh_{k}-t\\
1&i>Nh_{k}-t
\end{array}\right.
\ee
\bl
$\underset{N\rightarrow
\infty}{\lim} \d^{(N)}_{k}(T^{t})$ exists.
\el
\begin{proof} Suppose $N_{0}>>th_{k}$, let $\eta ^{(t)}_{k}=\frac{1}{th_{k}}\sum ^{th_{k}}_{i=1}d^{N,k}_{code}(i,t)$\\
\\
$\;$For any $n>N_{0}$, suppose $n=tp_{0}+n'$ ($0\leq n'<t$), then\\

$\d^{(n)}_{k}(T^{t})=\frac{\displaystyle 1}{\displaystyle nh_{k}}\sum\limits^{nh_{k}}_{i=1}d^{(n,k)}_{code}(i,t)$\\
$\qquad\qquad\qquad$ $=\frac{\displaystyle tp_{0}}{\displaystyle n}\frac{\displaystyle 1}{\displaystyle tp_{0}h_{k}}(p_{0}\sum\limits^{th_{k}}_{i=1}d^{(n,k)}_{code}(i,t)+\sum\limits^{nh_{k}}_{i=tp_{0}h_{k}+1}d^{(n,k)}_{code}(i,t))$ \\
\\
Thus $$\frac{tp_{0}}{n}\eta^{(t)}_{k}\leq
\d^{(n)}_{k}(T^{t})\leq\frac{tp_{0}}{n}(\eta^{(t)}_{k}+\frac{b}{p_{0}h_{k}})$$
Therefore, it is obvious that:
\be
\underset{n\rightarrow\infty}{\lim}\d^{n}_{k}(T^{t})=\eta^{t}_{k}
\ee
Done
\end{proof}
Now let $\d_{k}(T^{t})=\underset{n\rightarrow \infty}{\lim} \d^{(n)}_{k}(T^{t})=\eta^{(t)}_{k}$\\

\bl
$\left|\d_{k}-\d_{k+1}(T^{t})\right|\leq\frac{\displaystyle 2t}{\displaystyle h_{k+1}}$\\
\el

Therefore $\d_{k}(T^{t})$ converges. Define
$$\d(T^{t})=\underset{k \rightarrow \infty}{\lim}\d_{k}(T^{t})$$
This derives the following proposition:
\bp
$\d (T^{h_{k}}, Id)<\frac{\displaystyle 1}{\displaystyle 2^{k}}$
\ep

\bc
$(\mX ,\cB ,\mu , T)$ is rigid.
\ec

\section{Relative Prime Relation}
In this section, we set up the number theoretical relation of any positive integer with the height
of the $k$-stack $h_{k}$. The result is more general than what is needed in section 3.\\

The sequence of integer $h_{k}$ is also described
by induction:
\begin{align} \label{hk}
\nonumber & h_{0}=1;\\
& h_{k+1}=b^{k+1}h_{k}+1
\end{align}

where $b\in \mN, b\geq 2$.\\

The notation $n(m), n,m\in \mN $ is used for the integer $k$ such that $k\equiv n \mod m$, and $0\leq k <m$.\\
\\
Let $\mZ_{q}=\mZ /q\mZ$ and $\mZ^{*}_{q}=\{c+q\mZ\|(c,q)=1\}$
, ( $(m,n)$ is the notation for the largest common divisor of $m$
and $n$ ). Thus $\mZ_{q}$ is a finite commutative ring, and
$\mZ^{*}_{q}\subset \mZ_{q}$ is the
multiplicative group of multiply invertible elements of $\mZ_{q}$.\\
\\
Given $b\in \mZ_{q}^{*}, \, b\geq 2$, we define a map from $\mZ_{q}^{*}\times \mZ_{q}$ to
itself ($T: \mZ_{q}^{*}\times \mZ_{q} \to \mZ_{q}^{*}\times \mZ_{q}$), which is
also understood as a skew product of the rotation
on $\mZ_{q}^{*}$. That is:
\be
T(x,y)=T_{q,b}(x,y)=(bx,xy+1)
\ee

It is well understood that $(bx,xy+1)\in \mZ_{q}^{*}\times R$

\bl
Suppose $b,q\in \mN$, $(b,q)=1$. then there exists $n=n(q,b)\in\mN$ such that
$$T^{sn}=Id, for some s\in \mN$$.
\el
\begin{proof}
Without loss of generality, suppose $b\in \mZ_{q}^{*}, \, b\geq 2$. Since $\mZ_{q}^{*}\times \mZ_{q}$ is a finite set, $T$ is just a permutation
with finite order.
\end{proof}

\bl \label{intr1}
For any $k\in \mN$, we have
\be
T^{k}(1,0)=(b^{k}(q), h_{k-1}(q))\equiv (b^{k},
h_{k-1})\mod q,\, k\in\mN
\ee
\el
\begin{proof}
It is easy to see that
$$T(1,0)=(b(q),1(q))=(b(q),h_{0}(q))$$
$$T^{k}(1,0)=T(b^{k-1}(q),h_{k-2}(q))=(b^{k}(q),(b^{k-1}h_{k-2}+1)(q))$$
$$\quad\quad\quad\quad\quad =(b^{k}(q), h_{k-1}(q))$$
Therefore the lemma is proved by induction.
\end{proof}
\bl \label{intrt2}
Suppose $(b,q)=1$ then $h_{k}\equiv 0 \mod q$ for infinitely many $k\in \mN$.
\el
\begin{proof}
By Lemma \ref{intr1}, $T^{sn}(1,0)=(1,0)$.\\
Therefore we have
$$T^{sn}(1,0)=(b^{sn}(q), h_{sn1-1}(q)) \equiv (1,0),\mod q$$
$$h_{sn+1}\equiv 0(q), \quad s\in \mN , n=n(b,q)$$
\end{proof}
Now we can reach the following conclusion based on of Lemma \ref{intrt2}:
\bp \label{prime} Given $b\in \mN,\, b \geq 2$, and the sequence ${h_{k}}$ defined by \ref{hk}. Then for any $q\in \mN$, there exist infinitely many
$k\in\mN$ such that
$$h_{k+1}\equiv 1 \mod q$$
\ep
\begin{proof}
Suppose $b$ is divided by $q$, it is done.\\
So we only need to investigate the case of  $1\leq d=(b,q)<q$, $q=dq'$, $1<q'\leq q$, $(b, q')=1$.\\
By Lemma \ref{intrt2},  $h_{k}\equiv 0 \mod q'$ for infinitely many $k\in \mN$.\\
 On the other hand for $k$ large enough, $d$ divides $b^{k+1}$ thus
$$b^{k+1}h_{k}\equiv 0\mod q$$
$$h_{k+1}=b^{k+1}h_{k}+1\equiv 1\mod q$$
\end{proof}
\section{A family of totally rank one maps}
Now we revisit the generalized Del Junco-Rudolph's map $(\mX , \cB
, \mu , T)$. In this section $q$ is a
given integer greater than $1$, acting as the power index of the map.\\
\\
Since the heights of the stack structure satisfies
$h_{k+1}=b^{k+1}h_{k}+1$, Proposition \ref{prime} shows that for
infinitely many $n$, $h_{n}\equiv 1(\mod q)$ for
infinitely many n. We use the notion $\cN_{q}$ to denote the
set of those integers, that is $\cN_{q}=\{k|k\in \mN, h_{k}\equiv 1\mod q\}$.

\bd [Totally Rank One] \label{Total}

We say a finite measure-preserving system $(\mX, \b, T, \mu)$ is totally rank one, if all positive integer powers of $T$ are rank one.

\ed

\bl \label{uniq} Let the sequence ${h_{k}}$ be defined by \ref{hk}, for any $n\in \cN_{p}$, and $i\neq j, \, 0 \leq i,j <h_{n}$, we have:
$qi\neq qj(\mod h_{n}), i\neq j,
0\leq i,j<h_{n}$.
\el

\begin{proof}
Since $h_{n}$ is relatively prime with $q$.
\end{proof}
Now suppose $B^{*}_{N}$ to be the base set in the $N$th stack column. Then $B_{N}^{*}$ is corresponding to the
cylinder set $C_{N}$ with name $B_{N}$(the $N$-block in the recursive formula).\\
\\
We know that
\be \label{ineq2}\mu (T^{h_{N}}(B_{N}^{*})\Delta
(B_{N}^{*}))<\frac{1}{b^{N-1}}\mu (B_{N}^{*})
\ee

and
\be
\mu
(\overset{h_{N}-1}{\underset{i=0}{\cup}}(T^{i}B_{N}^{*}))>1-bh_{N}\mu
(\cS_{k})
\ee
By Proposition \ref{ineq1} we have\
\be \label{ineq3}  \mu
(\overset{h_{N}-1}{\underset{i=0}{\cup}}(T^{i}B_{N}^{*}))>1-b^{2}h_{N}/h_{N+1}>1-1/b^{N-1}
\ee

Equation \ref{ineq3} shows that the stack of disjoint union of
$\{T^{i}B^{*}_{N}\}$ is almost the whole space, except for a part
of measure no more than $\frac{\displaystyle  1}{\displaystyle
b^{N-1}}$. Equation \ref{ineq2} tells us the first return time of $B^{*}_{N}$
is $h_{N}$ except for a set of measure no more than
$\frac{\displaystyle 1}{\displaystyle b^{N-1}}$. Now suppose $N\in
\cN_{p}$. let $A^{*}_{N}=B^{*}_{N}\cap T ^{-h_{N}}(B^{*}_{N})\cap
T^{-2h_{N}}(B^{*}_{N})\cap\cdots\cap T^{-(q-1)h_{N}}(B^{*}_{N})$.
It is easy to see that
$$\mu (A^{*}_{N})\geq \mu (B^{*}_{N})-\sum^{q-1}_{i=1}\mu (B^{*}_{N}\Delta T^{-ih_{N}}(B_{N}^{*}))$$
$$\qquad\qquad\quad\quad\quad\quad>\mu (B^{*}_{N})-\frac{(q-1)(q-2)}{2}\mu (T^{h_{N}}(B^{*}_{N})\Delta(B^{*}_{N}))$$
By \ref{ineq2} we know that
\be
\mu (A^{*}_{N})>
(1-\frac{\displaystyle (q-1)(q-2)}{\displaystyle 2}\cdot \frac{\displaystyle 1}{\displaystyle b^{N-1}})\mu (B^{*}_{N})
\ee
Define the function $\tau :\{1,2, \cdots , h_{N}\}$ by $\tau
(i)=iq(h_{N})$. Since $N\in \cN_{p}$, $h_{N}$ is relatively prime
to $q$, therefore by Lemma \ref{uniq}, $\tau (i)$ is a
$h_{N}$-permutation. it is obvious to see that
$T^{qi}(A^{*}_{N})\subset T^{\tau (i)}(B_{N}^{*})$ $(0< i\leq h_{N})$, therefore, we have the following 3 claims:\\
\renewcommand{\labelitemi}{$\bullet$}
\begin{itemize}
\item [(I)] $\{T^{qi}(A^{*}_{N}), 0\leq i<h_{N}\}$ is a collection of pairwise disjoint sets;

\item [(II)]
\begin{align}\label{disjoint}\nonumber
\mu ((T^{q})^{h_{N}}(A^{*}_{N})\Delta
A^{*}_{N}) & <2\mu (A^{*}_{N}\Delta B^{*}_{N})+\mu (B^{*}_{N}\Delta
(T^{q})^{h_{N}}(B^{*}_{N}))\\
\nonumber  & \leq (q-1)(q-2) \frac{1}{b^{N-1}}\mu (B^{*}_{N})+\frac{q}{b^{N-1}}\mu (B^{*}_{N})\\
   & =\frac{(q^{2}-2q+2)}{b^{N-1}}\mu (B^{*}_{N})
\end{align}
\item [(III)] we see that (I) and (II), together with \ref{ineq3}, imply
\begin{align} \label{ineq5} \mu
\nonumber (\overset{h_{N}-1}{\underset{i=0}{\cup}}(T^{q})^{i}(A^{*}_{N}))& >(1-\frac{(q-1)(q-2)}{2})\frac{1}{b^{N-1}}\mu
(\overset{h_{N}-1}{\underset{i=0}{\cup}}T^{i}(B^{*}_{N}))\\
 & =(1-\frac{(q-1)(q-2)}{2}\frac{1}{b^{N-1}})h_{N}\mu (B^{*}_{N})
\end{align}
\end{itemize}

We know that
$C_{N}=\overset{h_{N}-1}{\underset{i=0}{\cup}}T^{i}(B^{*}_{N})$ is
the union of all the levels in the $N$th stack, so $\mX-C_{N}$ is
the remainder of the spacer set taken away the set of spacers used
in the first $N$ steps during the cutting and stacking process.
Thus $\mX-C_{N}=\cS _{N}$, $\mu (C_{N})=1-\mu
(\cS_{N})>1-\frac{\displaystyle 2}{\displaystyle h_{N+1}}$

\begin{align}
\mu (\overset{h_{N}-1}{\underset{i=0}{\cup}}(T^{q})^{i}(A^{*}_{N}))& >(1-\frac{(q-1)(q-2)}{2}\frac{1}{b^{N-1}})(1-\frac{2}{h_{N+1}})
\end{align}
\br Equation \ref{ineq5} shows that the measure of $(\overset{h_{N}-1}{\underset{i=0}{\cup}}(T^{q})^{i}(A^{*}_{N}))$, the $T^{q}$-stack with base $A^{*}_{N}$, is close to the full measure, since $B^{*}_{N}$ is the base of the $N$th $T$-stack and $h_{N}$ is the corresponding height.
\er
Now, (I), (II) and (III) tell us the following:
\bt \label{main} All the notations as above, every generalized Del Junco-Ruldoph's map $T$ is totally rank one.
\et

\br Though the rigidity of $T^{q}$ is
derived from b), it is a simple implication of the fact that all
powers of a rigid automorphism
on a standard Borel space are rigid.
\er

\bc In the weak closure of each general del
Junco-Rudolph's map, there is a dense $G_{\delta}$ subset of rank
one transformations.
\ec


\begin{thebibliography}{999}
\bibitem{BAX} J.R. Baxter, \emph{A class of ergodic transformations having simple spectrum}, Prodeeding of AMS, VOL 27, No. 2, 1971.
\bibitem{CHA} R. V. Chacon, \emph{Weakly mixing transformations which are not strongly mixing}, Proc. Amer. Math. Soc. 22 1969 559--562.
\bibitem{FER1} S. Ferenczi, \emph{Systems of finite rank}, Collaq. Math. 73 (1997), P. 35-65.
\bibitem{GLA} E. Glasner, \emph{Ergodic theory via joinings}, Mathematical Surveys and Monographs, 101. American Mathematical Society, Providence, RI, 2003.
\bibitem{GOO1} G.R. Goodson, \emph{Ergodic Dynamical Systems Conjugate to Their Composition Squares} Acta Math. Univ. Comenianae, Vol. LXXI, 2(2002), P. 201-210.
\bibitem{JUN1} A. del Junco, D. Ruldoph \emph{A rank-one, rigid, simple, prime map}, Ergodic Theory and Dynamical Systems, (1987), 7, P.229-247.
\bibitem{JUN2} A. del Junco, D. Ruldoph \emph{On Ergodic action whose self-joinings are graphs}, Ergodic Theory and Dynamical Systems, (1987), 7, P.531-557.
\bibitem{JUN3} A. del Junco, \emph{Disjointness of measure-presearving transformations}, self-joining and category. Ergodic Theory and Dynamical Systems I, Proceedings of Special Year, (1981) P. 81-89.
\bibitem{PRI1} A.A. Prikhod'ko, V.V. Ryzhikov, Several questions and hypotheses concerning the limit polynomials for chacon transformation, arXiv:1207.0614, 2012.
\bibitem{VEE3} W.A. Veech, \emph{The Metric Theory of Interval Exchange
Transformation I : Generic Spectral Properties}, American Journal of
Mathematics, 107(6):1331-1359,1984.
\bibitem{VIA} M. Viana, \emph{Ergodic Theory of Interval Exchange maps}, Rev. Mat. Complut., 19(2006), no. 1, 7-100.
\bibitem{WUphd} Y. Wu, \emph{Applications of Rauzy Induction on the generic ergodic theory of interval exchange transformations}, Doctor of Philosophy Thesis, Rice University Electronic Theses and Dissertations(2006).
\bibitem{WUetd} Y. Wu \emph{Whirly 3-interval exchange transformations}, Ergodic Theory and Dynamical Systems, available on CJO2015. doi:10.1017/etds.2015.63.
\bibitem{WU} Y. Wu, D. Li, D. Li, Y. Wang,  \emph{Totally Rank One Interval Exchange Transformations}, Arxiv, Cornell University, arXiv:1604.02638, (2016).

\end{thebibliography}
\end{document}